\documentclass[11pt]{article}

\usepackage{eurosym}
\usepackage{amsthm,amssymb,amsmath}
\usepackage{epsfig}
\usepackage{pgf,tikz}

\tikzstyle{black_V}=[fill=black, draw=black, shape=circle, scale=0.3]
\tikzstyle{none}=[fill=none, draw=none, shape=circle, scale=0.6]

\newtheorem{theorem}{Theorem}
\newtheorem{coro}[theorem]{Corollary}
\newtheorem{lemma}[theorem]{Lemma}
\newtheorem{pro}[theorem]{Proposition}
\newtheorem{example}[theorem]{Example}
\newtheorem{definition}[theorem]{Definition}
\newtheorem{rem}[theorem]{Remark}

\author{\small  J.C. VALENZUELA-TRIPODORO\textsuperscript{1}, M.A. MATEOS-CAMACHO\textsuperscript{2}, 
\\ \small M. CERA\textsuperscript{3}, R.M. CASABLANCA\textsuperscript{4}, 
M.P. \'ALVAREZ-RUIZ\textsuperscript{5}\\\small
\textsuperscript{1} \small Dpto. Matem\'aticas. ETSIA. Universidad de Cádiz. 11202 Algeciras. SPAIN\\
\textsuperscript{2} \small Escuela Internacional de Doctorado. Universidad de Sevilla. 41013  Sevilla. SPAIN\\
\textsuperscript{3} \small Dpto Matem\'atica Aplicada I. ETSIA. Universidad de Sevilla. 41013 Sevilla. SPAIN\\
\textsuperscript{4} \small Dpto Matem\'atica Aplicada II. ETSI. Universidad de Sevilla. 41092 Sevilla. SPAIN\\
\textsuperscript{5} \small Dpto. Estad\'istica e I.O. ETSIA. Universidad de Cádiz. 11202 Algeciras. SPAIN
}

\title{  $p$-Strong Roman Domination in Graphs}

\begin{document}

    \maketitle
    \begin{abstract}  
Domination in graphs is a widely studied field, where many different definitions have been introduced
in the last years to respond to different network requirements. This paper presents a new dominating
parameter based on the well-known strong Roman domination model. Given a positive integer $p$, we call
a $p$-strong Roman domination function ($p$-StRDF) in a graph $G$ to a function
$f:V(G)\rightarrow  \{0,1,2, \ldots , \left\lceil \frac{\Delta+p}{p} \right\rceil \}$ having the
property that if $f(v)=0$, then there is a vertex $u\in N(v)$ such that $f(u) \ge 1+ \left\lceil
\frac{ |B_0\cap N(u)|}{p} \right\rceil $, where $B_0$ is the set of vertices with label $0$. The
$p$-strong Roman domination number $\gamma_{StR}^p(G)$ is the minimum weight (sum of labels) of a
$p$-StRDF on $G$. We study the NP-completeness of the \emph{$p$-StRD}-problem, we also
provide general and tight upper and lower bounds depending on several classical invariants of the graph
and, finally, we determine the exact values for some families of graphs.
    \end{abstract}

\vspace{2ex}

\section{Introduction and notation}

Throughout this paper, $G$ is a simple graph with vertex set $V=V(G)$ and edge set $E=E(G)$. For a vertex 
$v\in V$, the \emph{open neighborhood} $N(v)$ is the set $\{u\in V(G):uv\in E(G)\}$ and the \emph{closed
neighborhood} of $v$ is the set $N[v]=N(v)\cup\{v\}$. The \emph{degree} of a vertex $v\in V$ is $d_{G}(v)=|N(v)|$.
The \emph{minimum} and \emph{maximum degree} of a graph $G$ are denoted by $\delta=\delta(G)$ and
$\Delta=\Delta(G)$, respectively. We denote by $P_{n}$ the \emph{path} of order $n$, $C_{n}$ for the 
\emph{cycle} of length $n$ and $\overline{K_{n}}$ for the edgeless graph with $n$ vertices.

A \textit{leaf} of $G$ is a vertex of degree one, while a \textit{support vertex} of $G$ is a vertex 
adjacent to a leaf. An \emph{$S$-external private neighbor} of a vertex $v\in S$ is a vertex 
$u\in V\setminus S$ adjacent to $v$ but no other vertex of $S$. The set of all $S$-external private neighbors 
of $v\in S$ is called the \emph{$S$-external private neighborhood} of $v$ and is denoted $\mathrm{epn}(v,S)$.

A \emph{tree} is a connected graph containing no cycles. A tree $T$ is called a \emph{double star} if it contains 
exactly two vertices that are not leaves. A double star with respectively $p$ and $q$ leaves attached at each 
support vertex is denoted by $S_{p,q}.$

Given two different graphs $G$ and $H$, let us denote by $G\vee H$ the graph obtained by adding to $G\cup H$ 
all possible edges joining a vertex in $G$ with a vertex in $H$.

A set $S\subseteq V$ in a graph $G$ is called a \textit{dominating set} if every vertex of $G$ is either 
in $S$ or adjacent to a vertex of $S.$ The \textit{domination number} $\gamma(G)$ equals the minimum 
cardinality of a dominating set in $G$.

The concept of Roman domination has arisen as a solution to a classic problem of military defensive 
strategy introduced by Stewart~\cite{Stewart} and by Revelle and Rosing~\cite{ReVelle}, having its 
origin in the time of Emperor Constantine I. At that moment, the Roman Empire had more conquered cities 
than legions for their defense in case of an attack. Defending a city was enough with a legion, which 
might be positioned in such a city or moved from another neighboring city. Then, Emperor Constantine I 
decreed a strategy based on two facts: first, any unprotected city should be able to be defended by a 
neighboring city and, second, no legion could come to defend an attacked neighboring city if such legion 
left unguarded its original location. The goal was to minimize the costs of settlement and mobilization 
of the legions while guaranteeing the possibility of defense of each position of the empire. For this, 
they could place up to two legions in each military settlement.

The first formal definition of Roman domination was introduced in $2004$ by Cockayne et al.~\cite{CDHH}, 
inspired by the works mentioned above. A function $f:V(G)\rightarrow\{0,1,2\}$ is a \textit{Roman 
dominating function} (RDF) on $G$ if every vertex $u\in V$ for which $f(u)=0$ is adjacent to at least one 
vertex $v$ for which $f(v)=2$. The weight of an RDF is the value $f(V(G))=\sum_{u\in V(G)}f(u).$ The 
\textit{Roman domination number} $\gamma_{R}(G)$ is the minimum weight of an RDF on $G$. Afterwards, the 
properties of this invariant have been extensively studied.

In recent years, other variations of Roman domination have been introduced, generally modifying the 
conditions in which the vertices are dominated, or adding some additional property to the classic version 
of the Roman domination. We highlight, for example: the {\it independent Roman domination}~\cite{aejs}, 
the {\it maximal Roman domination}~\cite{asstv}, the {\it weak Roman domination}~\cite{hh}, the {\it edge 
Roman domination}~\cite{rn}, the {\it total Roman domination}~\cite{ahsy}, the {\it signed Roman
domination}~\cite{ahlzs}, the {\it mixed Roman domination}~\cite{ahv}, the {\it Roman 
$\{2\}$-domination}~\cite{chha} or the work by Hsu et al.~\cite{hsu} on parallel algorithms 
for connected domination problems on interval and circular-arc graphs.

In all previous variants of Roman domination, it is assumed that a legion is enough to defend a position 
from individual attacks. However, there may be situations in which, even for individual attacks, this 
defensive strategy is insufficient. Other new variants are defined in~\cite{aacsv} and~\cite{bhh}, to 
contemplate other situations. This approach can have vast applications in service network modeling such 
as distribution, maintenance or provisioning.

The attack capacity is increasingly wider, giving rise to new situations, for instance, when the attacks 
occur simultaneously. In such cases, the previous defensive strategies are weak and insufficient. Some works 
try to solve this problem. Henning in~\cite{h}, based on weak Roman domination, provides a new version of 
the defense of the Roman Empire against multiple and sequential attacks. Recently, in~\cite{km}, the protection 
of a graph against sequential attacks on its vertices or edges is studied, by positioning mobile guards on 
vertices according to certain structures, for example, an eternal dominating set.

However, many real situations remain unresolved with these models, since nowadays, the attacks can be 
multiple and also simultaneous as, for instance, fires with several sources, synchronous natural disasters 
in different areas, joint attacks in cybernetics or security systems, etc. Several works address this approach. 
In $2009$, the \textit{Roman k-domination}, for $k\geq 1$, was defined~\cite{kv} to provide a reply against 
$k$ attacks in different vertices of a graph. A function $f:V(G)\rightarrow\{0,1,2\}$ is a 
\textit{Roman k-dominating function} on $G$ if every vertex $u\in V$ for which $f(u)=0$ is adjacent to at 
least $k$ vertices, $v_{1},\ldots,v_{k}$, with $f(v_i)=2$, for $i=1,\ldots,k$. This defensive strategy is 
dependable, as a defenseless vertex is covered by $k$ neighbors, but it can sometimes be excessive or 
unnecessary.

In $2017$, the \textit{strong Roman domination} was introduced~\cite{amsvy} as a reinforcement of the Roman
domination against multiple and simultaneous attacks, where legions are placed in strong vertices to defend
themselves and, at least, half of its unsafe neighbors. For a graph $G$ of order $n$ and maximum degree $\Delta$,
let $f:V(G) \rightarrow \{ 0,1, \ldots ,\left\lceil \frac{\Delta}{2}  \right\rceil +1 \}$ be a function that 
labels the vertices of $G$. Let $B_j=\{ v \in V : f(v) = j \}$ for $j= 0,1$ and let $B_2=V\smallsetminus 
(B_0 \cup B_1) =\{v\in V : f(v)\ge 2 \}$. Then, $f$ is a strong Roman dominating function (StRDF) on $G$, 
if every $v\in B_0$ has a neighbour $u$, such that $u\in B_2$ and $f(u)\ge 1+ \left\lceil 
\frac{|N(u)\cap B_0|}{2}  \right\rceil$. The minimum weight, $w(f)=f(V)= \sum_{u\in V} f(u)$, over all the 
strong Roman dominating functions for $G$, is called the {\it strong Roman domination number} of $G$ and we 
denote it by $\gamma_{StR}(G)$. An StRDF of minimum weight is called a $\gamma_{StR}(G)$-function. After 
this work, many relevant contributions on the strong Roman domination~\cite{cs,cs_tree, mnb,pajk,xw} have 
been provided.

Despite the particular interest of the strong Roman domination strategy and the further development of this 
kind of study, we may observe that the definition is certainly restrictive. In the StRD model, we consider
simultaneous attacks to unprotected neighbors of strong vertices under the condition that the stronger vertex 
may defend, at least, one-half of its neighbors.

In this paper, we introduce \emph{$p$-strong Roman domination}, a refined strategy of strong Roman domination 
that relaxes the definition, allowing for the development of less expensive defensive strategies.

\begin{definition}
Given a positive integer $p$, a function $f:V(G) \rightarrow \{ 0,1, \ldots ,\left\lceil \frac{\Delta+p}{p}
\right\rceil \}$ is a \emph{$p$-strong Roman dominating function} (\emph{$p$-StRDF}) if for every vertex 
$u\in B_0$ there is a vertex $v\in N(u)$ such that $f(v)\ge 1+ \left \lceil \frac{|N(v)\cap B_0|}{p}
\right\rceil$. The minimum weight of such a function is called the \emph{$p$-strong Roman domination number} 
of the graph and it is denoted by $\gamma_{StR}^p(G)$.
\end{definition}

In other words, the strong Roman domination model ensures that each strong vertex is capable of defending at 
least half of its undefended neighbors without leaving its own location unprotected. This means that it has 
one unit to protect every group of two weak neighbors. In this case, the $p$-StRD model ensures that each 
strong position has at least one legion to defend each group of $p$ undefended neighbors.

\begin{figure}
\centering
\begin{tikzpicture}[scale=0.32]
		\node [style={black_V}] (0) at (-10, 0) {};
		\node [style={black_V}] (1) at (-12, 3) {};
		\node [style={black_V}] (2) at (-10, 4) {};
		\node [style={black_V}] (3) at (-8, 3) {};
		\node [style={black_V}] (4) at (-13, 1) {};
		\node [style={black_V}] (5) at (-13, -1) {};
		\node [style={black_V}] (6) at (-12, -3) {};
		\node [style={black_V}] (7) at (-10, -4) {};
		\node [style={black_V}] (8) at (-8, -3) {};
		\node [style={black_V}] (9) at (-7, -1) {};
		\node [style={black_V}] (10) at (-7, 1) {};
		\node [style=none] (33) at (-10.30, 1) {\large $2$};
		\node [style=none] (34) at (-10.5, 4.25) {\large $0$};
		\node [style=none] (35) at (-8.5, 3.5) {\large $0$};
		\node [style=none] (36) at (-7.5, 1.25) {\large $0$};
		\node [style=none] (37) at (-7.25, -0.5) {\large $0$};
		\node [style=none] (38) at (-8.5, -2.75) {\large $0$};
		\node [style=none] (39) at (-10.5, -3.75) {\large $0$};
		\node [style=none] (40) at (-12.5, -2.75) {\large $0$};
		\node [style=none] (41) at (-13.5, -0.75) {\large $0$};
		\node [style=none] (42) at (-13.5, 1.25) {\large $0$};
		\node [style=none] (43) at (-12.5, 3.25) {\large $0$};
		\node [style={black_V}] (44) at (-1, 0) {};
		\node [style={black_V}] (45) at (-3, 3) {};
		\node [style={black_V}] (46) at (-1, 4) {};
		\node [style={black_V}] (47) at (1, 3) {};
		\node [style={black_V}] (48) at (-4, 1) {};
		\node [style={black_V}] (49) at (-4, -1) {};
		\node [style={black_V}] (50) at (-3, -3) {};
		\node [style={black_V}] (51) at (-1, -4) {};
		\node [style={black_V}] (52) at (1, -3) {};
		\node [style={black_V}] (53) at (2, -1) {};
		\node [style={black_V}] (54) at (2, 1) {};
		\node [style=none] (55) at (-1.25, 1) {\large $6$};
		\node [style=none] (56) at (-1.5, 4.25) {\large $0$};
		\node [style=none] (57) at (0.5, 3.5) {\large $0$};
		\node [style=none] (58) at (1.5, 1.25) {\large $0$};
		\node [style=none] (59) at (1.75, -0.5) {\large $0$};
		\node [style=none] (60) at (0.5, -2.75) {\large $0$};
		\node [style=none] (61) at (-1.5, -3.75) {\large $0$};
		\node [style=none] (62) at (-3.5, -2.75) {\large $0$};
		\node [style=none] (63) at (-4.5, -0.75) {\large $0$};
		\node [style=none] (64) at (-4.5, 1.25) {\large $0$};
		\node [style=none] (65) at (-3.5, 3.25) {\large $0$};
		\node [style={black_V}] (66) at (8, 0) {};
		\node [style={black_V}] (67) at (6, 3) {};
		\node [style={black_V}] (68) at (8, 4) {};
		\node [style={black_V}] (69) at (10, 3) {};
		\node [style={black_V}] (70) at (5, 1) {};
		\node [style={black_V}] (71) at (5, -1) {};
		\node [style={black_V}] (72) at (6, -3) {};
		\node [style={black_V}] (73) at (8, -4) {};
		\node [style={black_V}] (74) at (10, -3) {};
		\node [style={black_V}] (75) at (11, -1) {};
		\node [style={black_V}] (76) at (11, 1) {};
		\node [style=none] (77) at (7.75, 1) {\large $4$};
		\node [style=none] (78) at (7.5, 4.25) {\large $0$};
		\node [style=none] (79) at (9.5, 3.5) {\large $0$};
		\node [style=none] (80) at (10.5, 1.25) {\large $0$};
		\node [style=none] (81) at (10.75, -0.5) {\large $0$};
		\node [style=none] (82) at (9.5, -2.75) {\large $0$};
		\node [style=none] (83) at (7.5, -3.75) {\large $0$};
		\node [style=none] (84) at (5.5, -2.75) {\large $0$};
		\node [style=none] (85) at (4.5, -0.75) {\large $0$};
		\node [style=none] (86) at (4.5, 1.25) {\large $0$};
		\node [style=none] (87) at (5.5, 3.25) {\large $0$};
		\node [style=none] (88) at (-10, -5.5) {\Large $(a)$};
		\node [style=none] (89) at (-1, -5.5) {\Large $(b)$};
		\node [style=none] (90) at (8, -5.5) {\Large $(c)$};

		\draw (10) to (0);
		\draw (0) to (9);
		\draw (0) to (8);
		\draw (0) to (7);
		\draw (0) to (6);
		\draw (0) to (5);
		\draw (0) to (4);
		\draw (0) to (1);
		\draw (2) to (0);
		\draw (0) to (3);
		\draw (54) to (44);
		\draw (44) to (53);
		\draw (44) to (52);
		\draw (44) to (51);
		\draw (44) to (50);
		\draw (44) to (49);
		\draw (44) to (48);
		\draw (44) to (45);
		\draw (46) to (44);
		\draw (44) to (47);
		\draw (76) to (66);
		\draw (66) to (75);
		\draw (66) to (74);
		\draw (66) to (73);
		\draw (66) to (72);
		\draw (66) to (71);
		\draw (66) to (70);
		\draw (66) to (67);
		\draw (68) to (66);
		\draw (66) to (69);
\end{tikzpicture}
\caption{For a star graph: (a) a RDF, (b) an StRDF, and (c) a $4$-StRDF.}
\end{figure}
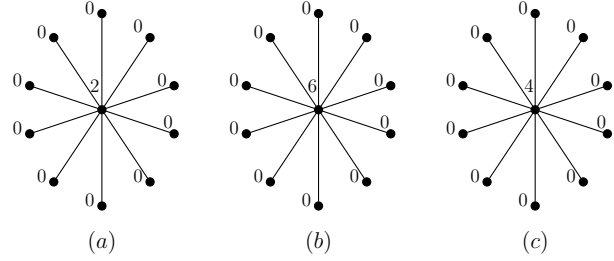 \label{fig:example}


\section{Complexity results}


This section aims to establish the NP-completeness of the $p$-StRD problem for bipartite 
and chordal graphs. The following decision problem is associated with the optimization problem 
of calculating the $p$-StRD number of a given graph.

\bigskip 

\noindent{\bf \emph{pStRD-Number Problem}}\newline
\textbf{Instance}: Graph $G =(V,E)$ and a positive integer $r$.\newline
\textbf{Question}: Does $G$ have a $p$-StRD function $f$ with $f(V)\le r$?

\bigskip

We make use of the \emph{Exact Cover by 3-Sets (X3C) problem} (see \cite{GJ1979}) to demonstrate 
that \emph{pStRD-Number Problem} is NP-complete. Namely, an instance of X3C is the following

\bigskip

\noindent{\bf EXACT 3-Cover (X3C) Problem}\newline
\textbf{Instance}: A collection $C$ of $3$-element subsets of a finite set $X$ with $|X|=3q.$\newline
\textbf{Question}: Does $X$ have an \emph{exact cover} in $C,$ that is, a subcollection $C'\subseteq C$ 
that contains every element of $X$ in exactly one member?

{ 
{\bf Example. } Let \( q=2 \text{ and } X = \{x_1, x_2, \dots, x_6\} \) be a set of \( 3q \) literals, and let \( C=\{(x_1,x_2,x_3),\) 
\((x_1,x_2,x_4),\) \( (x_1,x_5,x_6), \) \( (x_2,x_3,x_4), \) \( (x_3,x_5,x_6)\} \) be a collection of clauses 
\emph{(}subsets of cardinality 3\emph{)} of \( X \). Clearly, \( C' = \{(x_1, x_2, x_4), (x_3, x_5, x_6)\} \)
\(\subseteq C\) is an exact cover of \( C, \) because each and every element of $X$ belongs to exactly one clause
in $C'$. Note that $|C'|=q=2.$
}

The main result of this section is presented in the following theorem.

\begin{theorem}\label{complexity_bipartite}  
The \emph{$p$-StRD} number problem is NP-complete, even when restricted to bipartite or chordal graphs.
\end{theorem}

\noindent \textbf{Proof.}

First, {  this problem belongs to the class of NP problems since we could verify, in polynomial 
time concerning the size $n$, whether a given possible solution is indeed a solution or not.}

Next, we prove that \emph{pStRD}-Number Problem is NP-complete for bipartite graphs {  by constructing 
a polynomial-time transformation from X3C problem,  which is a well-known NP-problem (see \cite{GJ1979}).}

Assume that $I=(X,C)$ is an arbitrary instance of X3C, with $X=\{x_{1},x_{2},\ldots ,x_{3q}\}$ and
$C=\{C_{1},C_{2},\ldots ,C_{t}\}$. 

{  The key of this proof is that we construct a bipartite graph $B(I)$ starting from $I$} 
and provide a positive integer $r$ such that $I$ contains an \emph{exact cover} by $3$-sets if and only if $B(I)$ 
has a $p$-StRD function $f$ having weight $w(f)\le r=2q+3t$

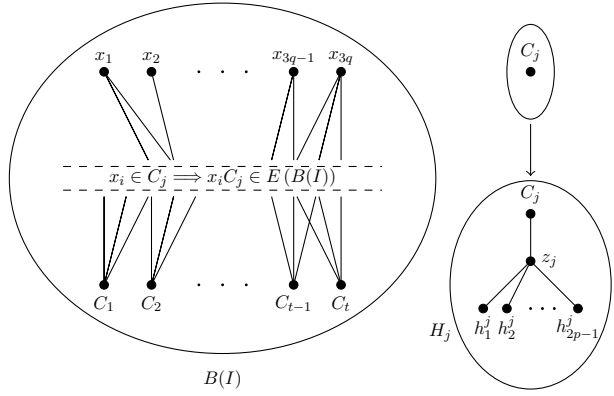
\begin{figure}[h]
\centering
\begin{tikzpicture}[scale=0.63, transform shape]
\node [draw, shape=circle,fill=black,scale=0.5] (x1) at  (1,0) {};
\node at (1,0.3) {$x_1$};
\node [draw, shape=circle,fill=black,scale=0.5] (x2) at  (2,0) {};
\node at (2,0.3) {$x_2$};
\node [draw, shape=circle,fill=black,scale=0.1] () at  (3,0) {};
\node [draw, shape=circle,fill=black,scale=0.1] () at  (3.5,0) {};
\node [draw, shape=circle,fill=black,scale=0.1] () at  (4,0) {};
\node [draw, shape=circle,fill=black,scale=0.5] (x3) at  (5,0) {};
\node at (5,0.3) {$x_{3q-1}$};
\node [draw, shape=circle,fill=black,scale=0.5] (x4) at  (6,0) {};
\node at (6,0.3) {$x_{3q}$};
\node  (x11) at  (2,-2) {};
\node  (x12) at  (2.5,-2) {};
\node  (x13) at  (4.5,-2) {};
\node  (x14) at  (5,-2) {};
\node  (x15) at  (5.5,-2) {};
\node  (x16) at  (6,-2) {};
\draw (x1)--(x11)--(x1)--(x12)--(x2);
\draw (x3)--(x13)--(x3)--(x14)--(x4)--(x15)--(x4)--(x16);
\node  (x00) at  (0,-2) {};
\node (x01) at  (7,-2) {};
\draw[dashed] (x00)--(x01);
\node  (x02) at  (0,-2.5) {};
\node (x03) at  (7,-2.5) {};
\draw[dashed] (x02)--(x03);
\node at (3.5,-2.25){$x_i\in C_j\Longrightarrow x_iC_j\in E\left(B(I)\right)$};
\node [draw, shape=circle,fill=black,scale=0.5] (c1) at  (1,-4.5) {};
\node at (1,-4.9) {$C_1$};
\node [draw, shape=circle,fill=black,scale=0.5] (c2) at  (2,-4.5) {};
\node at (2,-4.9) {$C_2$};
\node [draw, shape=circle,fill=black,scale=0.1] () at  (3,-4.5) {};
\node [draw, shape=circle,fill=black,scale=0.1] () at  (3.5,-4.5) {};
\node [draw, shape=circle,fill=black,scale=0.1] () at  (4,-4.5) {};
\node [draw, shape=circle,fill=black,scale=0.5] (c3) at  (5,-4.5) {};
\node at (5,-4.9) {$C_{t-1}$};
\node [draw, shape=circle,fill=black,scale=0.5] (c4) at  (6,-4.5) {};
\node at (6,-4.9) {$C_t$};
\node  (c11) at  (1,-2.5) {};
\node  (c12) at  (1.5,-2.5) {};
\node  (c13) at  (2,-2.5) {};
\node  (c14) at  (2.5,-2.5) {};
\node  (c15) at  (3,-2.5) {};
\node  (c16) at  (4.5,-2.5) {};
\node  (c17) at  (5,-2.5) {};
\node  (c18) at  (5.5,-2.5) {};
\node  (c19) at  (6,-2.5) {};
\draw (c1)--(c11)--(c1)--(c12)--(c1)--(c13)--(c2)--(c14)--(c2)--(c15);
\draw (c16)--(c3)--(c17)--(c4)--(c18)--(c3);
\draw (c4)--(c19);

\draw (3.5,-2.25) ellipse (4.5cm and 3.7cm);
\node at (3.5,-6.5) {$B(I)$};
%
%

\node [draw, shape=circle,fill=black,scale=0.5] (cj) at  (10,0) {};
\node at (10,0.4) {$C_j$};
\draw (10,0) ellipse (0.5cm and 1cm);

\node [draw, shape=circle,fill=black,scale=0.5] (cj1) at  (10,-3) {};
\node at (10,-2.6) {$C_j$};
\node [draw, shape=circle,fill=black,scale=0.5] (zj) at  (10,-4) {};
\node at (10.4,-4) {$z_j$};
\draw (cj1)--(zj);
\node [draw, shape=circle,fill=black,scale=0.5] (y11) at  (9,-5) {};
\node at (9,-5.4) {$h_1^j$};
\node [draw, shape=circle,fill=black,scale=0.5] (y12) at  (9.5,-5) {};
\node at (9.5,-5.4) {$h_2^j$};

\node [draw, shape=circle,fill=black,scale=0.1] () at  (10,-5) {};
\node [draw, shape=circle,fill=black,scale=0.1] () at  (10.25,-5) {};
\node [draw, shape=circle,fill=black,scale=0.1] () at  (10.5,-5) {};

\node [draw, shape=circle,fill=black,scale=0.5] (y14) at  (11,-5) {};
\node at (11,-5.4) {$h_{2p-1}^j$};
\draw (zj)--(y11);
\draw (zj)--(y12);
\draw (zj)--(y14);

\draw (10,-4.5) ellipse (1.7cm and 2.2cm);
\draw[->] (10,-1.1)--(10,-2.2);
\node at (8.1,-5.5) {$H_j$};

\end{tikzpicture}
\caption{Constructing $B(I)$ and $C_j$-gadgets.}
\label{fig-1}
\end{figure}

Let $B(I)$ the bipartite graph with classes $X=\{ x_i: 1\leq i \le 3q\}$ and $C=\{C_j:1\le j\le t\}$ 
being $x_i$ adjacent to {$C_j$} if $x_i\in C_j$. Next, the vertices $C_j$ of the bipartite graph will be 
swapped out for an appropriate gadget.

For each $C_j\in C,$ let $H_j$ be the star with central vertex $z_j$ and leaves $\{C_j,h_1^j, \ldots , 
h_{2p-1}^j \}$. We replace $C_j$ with $H_j$ by identifying $C_j$ with the corresponding leaf $C_j$ of $H_j.$ 
The new graph, denoted by $B(I),$ is bipartite with classes
$\{x_i: 1\leq i \leq 3q\}\cup\{z_j: 1\leq j\leq t\}$ and
$\{C_j\ : 1\leq j\leq t\}\cup\{h_l^j\ : 1\leq j\leq t,\ 1\leq l\leq 2p-1\}\}$. We set $r=2q+3t.$
Of course, we can construct $B(I)$ in polynomial time on the size of the given instance.

By assuming that $C'$ is an exact cover for $X$ in $C$, we define the following function over 
$V(B(I))$:

{\small
$$ f(v)=\left\{
\begin{array}{ccl} 
0 & \text{if} & v \in \{x_i: 1\leq i \leq 3q\}\\[0.25em]
0 & \text{if} & v \in \{h_l^j: 1\leq j\leq t,\ 1\leq l\leq 2p-1\} \\[0.25em]
0 & \text{if} & v \in \{C_j: C_j\not\in C', 1\leq j\leq t\}\\[0.25em]
2 & \text{if} & v \in \{C_j: C_j\in C', 1\leq j\leq t\}\\[0.25em]
3 & \text{if} & v \in \{z_j: 1\leq j\leq t\}
\end{array}
\right.
$$
}

Clearly, $C'$ has cardinality equal to $q$, {  because $C'$ covers $C$,
$|C|=3q,$ and each clause in $C'$ has cardinality equal to $3$. Besides,} 
$f\left(\{C_j\ : C_j\in C', 1\leq j\leq t\}\right)=2q.$ 

Since $C'$ is an exact cover 
for $X$ in $C,$ we know that for all $1\leq i \leq 3q$ there exists $C_j\in C'$ with 
$x_i\in C_j.$ Therefore
$f(C_j)=2 \ge 1 + \left \lceil \frac{|N(C_j)\cap B_0|}{p}\right\rceil =
1 + \left \lceil\frac{3}{p} \right \rceil=2$ for $p\ge 3$. Furthermore, for each
$v\in \{C_j : C_j\not\in C', 1\leq j\leq t\}\cup \{h_l^j : 1\leq j\leq t,\ 1\leq l\leq 2p-1\},$ 
there exists $k\in \{1\ldots ,t\}$ with $z_k \in N(v)$ and
$f(z_k)=3\ge 1+ \left \lceil \frac{|N(z_k)\cap B_0|}{p}  \right\rceil=1+2.$ 

So, $f$ is a $p$-StRD function with $f(V(B(I)))=2q+3t=r.$

On the other hand, suppose now that there exists a $p$-StRD function $f$ with $f(V(B(I)))\le r.$
Without loss of generality, let 
$f$ be one of those functions that assigns as much label value as feasible to 
the set $\{z_j : 1\leq j\leq t\}$.

Under these conditions, we may readily verify that $f(z_j)=3$ and $f(h_l^j)=0$ for all
$1\leq j \leq t$ and $1\leq l \le 2p-1$. On the contrary, if there exists $j\in\{1,\ldots,t\}$ 
and $l\in \{1,\ldots,2p-1\}$ with $f(h_l^j)\geq 1$ then either $f(h_l^j)\ge 1$ for all $1\le l \le 2p-1$  
or $f(z_j)\ge 2$ because not all $f(h_l^j)\ge 1.$ Anyhow, we can define $f^*$ such that $f^*(z_j)=3$ 
and $f^*(h_l^j)=0.$  

Hence, it follows that { 
{\small%
$$ f\left(\{z_j\ : 1\leq j\leq t\}) \cup \{h_l^j: 1\leq j\leq t,\ 1\leq l\leq 2p-1\} \right)$$}%
is equal to $3t,$ and then\newline 
$f\left(\{c_j: 1\leq j\leq t\}) \cup \{x_i: 1\leq i \leq 3q\} \right)\le r-3t.$
}

Let us denote by $A_i = \{x \in X : f(x)=i\}  \text{ and }  a_i=|A_i|,$ for $i=0,1;$ 
$A_2 = \{x \in X : f(x)\ge 2\}  \text{ and }  a_2=|A_2|;
D_i = \{C_j \in C : f(C_j)=i\}  \text{ and }  d_i=|D_i|,$ for $i=0,1;$
$D_2 = \{C_j \in C : f(C_j)\ge 2\}  \text{ and }  d_2=|D_2|.$

The following equalities, $a_0+a_1+a_2=3q$ and $d_0+d_1+d_2=t$, are an immediate consequence of
this notation. Note that $a_1+2 a_2+ d_1 +2 d_2 
\le a_1+ f(A_2)+d_1+ f(D_2) \le 2q$ because $f(V)\le r.$ Additionally, for all
$x\in A_0$ there exists $C_j \in D_2$ with $x\in N(C_j)$ and then $d_2\ge \frac{a_0}{3}.$
We have that $ 2q\ge  a_1+2 a_2+ d_1 +2 d_2 =3q-a_0+a_2+ d_1 +2 d_2 \ge 3q-3d_2+a_2+ d_1 +2 d_2$
and therefore $d_2-(a_2+d_1)\ge q.$ So, $d_2\ge q,$ and since $a_1+2 a_2+ d_1 +2 d_2 \le 2q ,$ it
is verified that $d_2=q$ and $f(C_j)=2$ { for all $C_j\in D_2$ } and, finally { $a_1= a_2=d_1=0.$}
As a result, we get that $d_1=0$ and $d_0=t-q.$ Since $|X|=3q=3|D_2|,$ we have that the set 
$C'=\{C_j\ :  C_j \in D_2\}$ is a subcollection $C'\subseteq C$ that contains every element of $X$ 
in exactly one member of $C',$ and the result follows for bipartite graphs.

By adding all the edges between the vertices $C_j$'s, we obtain a chordal graph. Consequently, by using a similar proof to the one developed to arrive at the previous result, we can derive the one
for chordal graphs.

$\hfill \Box $

\section{General bounds}

This section is dedicated to presenting general and different bounds for the $p$-strong Roman domination 
number in graphs, which is a natural step after checking the NP-completeness of the $p$-StRD Roman 
domination problem in the previous section. 

First of all, let us see which values of parameter $p$ have to be considered.

Let $G$ be any graph of order $n$. Let $p$ be a positive integer and let $f$ be a $p$-StRD function having
minimum weight in the graph $G$. Let us see that $3\le p \le \Delta-1$ have to be assumed.

If $p=1$, then
$$ \begin{array}{rcl}
	w(f)  & =     &\displaystyle
			     \sum_{v\in B_1\cup B_2} f(v) = \sum_{v\in B_1} f(v) +\sum_{v\in B_2} f(v) \\[1em]
		&  \ge & \displaystyle
			     |B_1| + \sum_{v\in B_2} \left(1+ \left \lceil \frac{|N(v)\cap B_0|}{p}  \right\rceil\right) \\[1em]
		&   =   &  |B_1| + |B_2|+\displaystyle \sum_{v\in B_2}  |N(v)\cap B_0|. 
\end{array}	
$$

Since $f$ is a $p$-StRD function of minimum weight, each $v\in B_2$ must have a private neighbor in $B_0$
and therefore $w(f) \ge |B_1|+|B_2|+|B_0|=|V(G)|=n.$ Hence, for $p=1$, the function $f(u)=1$ for all $u\in V(G)$ 
is a $1$-strong Roman domination function of minimum weight and $\gamma_{StR}^p(G)=n.$

If $p=2$, taking into account the definition of a $p$-StRD function, we may derive that the strategy of 
$2$-strong Roman domination is just the same as the one of the strong Roman domination model.

Finally, if $p\ge \Delta$, then $\left\lceil \frac{\Delta}{p}  \right\rceil = 1$. Hence we have that 
$f:V(G) \rightarrow \{ 0,1, 2 \}$ and the condition $f(v)\ge 1+ \left \lceil \frac{|N(v)\cap B_0|}{p}  
\right\rceil$ is equivalent to $f(v)=2$, because
$$ 1<1+ \left \lceil \frac{|N(v)\cap B_0|}{p}  \right\rceil  \le  1+\left\lceil \frac{\Delta}{p} \right\rceil = 2$$
Therefore, $p$-strong Roman domination corresponds to the original Roman domination model when $p\ge \Delta$.

Summing up, the $p$-StRD model is trivial for $p=1$, matches the StRD strategy for $p=2$ and coincides with 
the original Roman domination problem for all $p\ge \Delta$. Therefore, from now on, we will only consider 
$3\le p \le \Delta-1$. Observe that the latter implies that $\Delta \ge 4.$


As an immediate consequence of the definition, we can point out the following remark.
\begin{rem}
Let $G$ be a connected graph having maximum degree $\Delta \ge 4.$ Let $p,q$ be positive integers 
such that $3\le p \le q \le \Delta -1$. Then, $$\gamma_{StR}^p(G) \ge \gamma_{StR}^q(G).$$
\end{rem}
Of course, it is not difficult to relate our new parameter to some of the most well-known parameters in domination. 
As an initial bound for the $p$-StRD number, we prove the following result.

\begin{rem} \label{otras}
Let $G$ be a connected graph having maximum degree $\Delta \ge 4.$ Let $p$ be a positive 
integer such that $3\le p \le \Delta -1$. Then,
$$ \gamma_R(G) \le \gamma_{StR}^p(G) \le \left( \left\lceil \frac{\Delta}{p}  \right\rceil + 1 \right) \gamma(G).$$
\end{rem}

\noindent  {\bf Proof.} For the lower bound, let $f=(B_0,B_1,B_2)$ be {  any $\gamma_{StR}^p(G)$-function} on $G$ and define the function
$g:V(G)\rightarrow\{0,1,2\}$, such that $g(u)=2$ whenever $u\in B_2$ and $g(u)=f(u)$ otherwise. Hence, 
$\gamma_R(G) = |V_1|+2|V_2| = |B_1|+2|B_2| \le \gamma_{StR}^p(G)$. On the other hand, let $D$ be a dominating set 
and let $f:V(G) \rightarrow \{ 0,1, \ldots ,\left\lceil \frac{\Delta}{p}  \right\rceil +1 \}$ be the function 
defined as follows $f(u)=\left\lceil \frac{\Delta}{p}  \right\rceil + 1$ for all $u \in D$ and $f(u)=0$ 
otherwise. The function {  $f$ is a $p$-StRD function}, which leads us to the upper bound. 

$\hfill \Box $

Next, we prove an upper bound that only depends on the order and the maximum degree of the graph.

\begin{pro}\label{up1}
Let $G$ be a graph with order $n$ and maximum degree $\Delta \ge 4$. Let $p$ be a positive integer 
such that $3\le p \le \Delta -1$. Then
$$ \gamma_{StR}^p(G) \le n-\Delta+\left\lceil \frac{\Delta}{p}  \right\rceil $$
\end{pro}

\noindent {\bf Proof.} Let $u$ be a vertex with degree {  $\Delta.$} Let us define the function {  $f:V(G) \rightarrow \{ 0,1, \ldots ,\left\lceil \frac{\Delta}{p}
\right\rceil + 1\}$ } as follows:
$f(u)=\left\lceil \frac{\Delta}{p}  \right\rceil + 1$, $f(v)=0$ for all $v\in N(u)$ and $f(v)=1$ otherwise. {  Taking into account that $B_2=\{u\}, B_0=N(u)$ and $B_1=V \setminus N[u],$ then $f$ is a $p$-StRD function and therefore}
$$ \begin{array}{rcl}
\gamma_{StR}^p(G) &\le &  w(f)\\ 
 &=&\displaystyle |B_1| + \sum_{x\in B_2} f(x)\\ [1.5em] 
 &=& |V \setminus N[u]| +f(u) \\ 
 &=& (n-\Delta-1) + \left\lceil \frac{\Delta}{p} \right\rceil+1\\ 
 &=& n-\Delta+\left\lceil \frac{\Delta}{p}  \right\rceil .
\end{array}$$
 \hfill $\Box$

It is worth noting that the upper bound given by Proposition~\ref{up1} is sharp, for example, for every 
star $K_{1,n-1}$ with $3\le p\le n-2.$

\begin{coro} Let $G$ be a graph with order $n$ and maximum degree $\Delta \ge 4$. Let $p$ be a positive 
integer such that $3\le p \le \Delta -1$. Then
$$ \gamma_{StR}^p(G) \le n-2 $$
\end{coro}

\noindent {\bf Proof.} By applying Proposition~\ref{up1}, we have that 
$$ 
\begin{array}{rcl}
\gamma_{StR}^p(G) &\le& n-\Delta+\left\lceil \frac{\Delta}{p}  \right\rceil
\le  n-\Delta+\frac{\Delta}{p}+1\\[1em] 
&\le& n+1-\frac{2\Delta}{3} \le n+1 -\frac{8}{3} \le n-2.
\end{array}$$
$ \hfill \Box$

Our next result concerns improving the previous bound for $r$-regular graphs.

\begin{pro} \label{regulares}
Let $3\le p \le \Delta -1$ be a positive integer and let $G$ be a $r$-regular graph, with $r \ge p+1$ 
and girth $g\ge 5.$ Then
$$\displaystyle \gamma_{StR}^p (G) \le n-r^2+ \left( \left\lceil \frac{r-1}{p}  \right\rceil  + 1\right)r$$
\end{pro}

\noindent {\bf Proof.}
Consider any vertex $u\in V(G)$. Since $G$ is an $r$-regular graph, we have that $N(u)$ is a set of 
$r$ vertices, say $N(u)=\{w_1,\ldots ,w_{r}\}.$ {  Moreover, each one of the sets $N(w_j)-u$, 
for $j=1 \ldots r$, is formed by $r-1$ different vertices, say $N(w_j)-u=\{ z^j_1,\ldots ,z^j_{r-1} \}.$ 
Note that, due to the girth of $G$, $N(u)$ is an independent set; each set $N(w_j)-u$ is also an independent set; 
and they are disjoint set of vertices of $V(G).$
}

Let us define a function $f$ 
as follows: $f(u)=1$, $f(w_j)= 1 + \left\lceil \frac{r-1}{p}  \right\rceil$ for all $1\le j\le r$, $f(z^j_k)=0$ 
for all $1\le k \le r-1$, and $f(v)=1$ for any non yet labelled vertex $v$, if any. {  As we have observed before, the 
vertices $u, w_j, z^j_k$ are all different, since the girth of $G$ is at least $5$. Finally, any vertex 
$z^j_k$ is dominated by a vertex $w_j$ labelled with a label equal to} $1 + \left\lceil \frac{r-1}{p}  \right\rceil =
1+\left\lceil \frac{|N(w_j)\cap B_0|}{p} \right\rceil$, therefore, the defined function $f=(B_0,B_1,B_2)$ is a $p$-StRD
function on $G$ and it holds
$$ 
\begin{array}{rcl}
\gamma_{StR}^p (G) & \le & w(f) = \displaystyle |B_1| + \sum_{x\in B_2} f(x) \\[1em]
& \leq &  1 + \left[ n - (1+r+(r-1)r) \right ]\\[1em] 
& &+ \left( 1 + \left\lceil \frac{r-1}{p}  \right\rceil \right) r \\[1em]
& = &  n-r^2+ \left( \left\lceil \frac{r-1}{p}  \right\rceil  + 1\right)r.
\end{array}
$$
$\hfill \Box$

\begin{coro} \label{regularesco}
Let $3\le p \le \Delta -1$ be a positive integer and let $G$ be a $(p+1)$-regular graph with girth $g \ge 5.$ Then
$$ \gamma_{StR}^p (G) \le n-p^2+1$$
\end{coro}

To check the tightness of this upper bound for regular graphs, we first prove a technical result which 
will be useful later.

\begin{lemma} \label{lowB0} Let $G$ be a graph with order $n$ and maximum degree $\Delta \ge 4$. Let $p$ 
be a positive integer such that $3\le p \le \Delta -1$. Let $f=(B_0,B_1,B_2)$ be any $p$-StRD function 
on $G$.  Hence
$$ \gamma_{StR}^p(G) \ge  n+ \left\lceil \frac{1-p}{p} |B_0| \right\rceil .$$
\end{lemma}

\noindent {\bf Proof.} Observe that each vertex in $B_1\cup B_2$ adds one unit, by itself, to the weight 
of $f$. In addition, since every vertex in $B_0$ has, at least, a neighbor in $B_2$, each vertex in $B_0$ 
adds, at least, $\frac{1}{p}$ units to the weight of $f$. Therefore
$$ 
\begin{array}{rcl}
\gamma_{StR}^p(G)&=&w(f) \ge |B_1|+ |B_2| + \left\lceil  \frac{1}{p}  \right\rceil |B_0|\\[1em]
&=& n-|B_0|+ \left\lceil  \frac{|B_0|}{p}  \right\rceil 
\ge n+ \left\lceil \frac{1-p}{p} |B_0| \right\rceil
\end{array}
$$
$\hfill \Box$

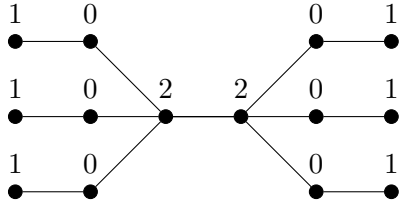
\begin{figure}[htb]
\center
\begin{tikzpicture}[scale=0.5]
\node [draw, shape=circle,fill=black,scale=0.5] (a1) at  (0,4) {};
\node [draw, shape=circle,fill=black,scale=0.5] (a2) at  (0,2) {};
\node [draw, shape=circle,fill=black,scale=0.5] (a3) at  (0,0) {};
\node [draw, shape=circle,fill=black,scale=0.5] (b1) at  (2,4) {};
\node [draw, shape=circle,fill=black,scale=0.5] (b2) at  (2,2) {};
\node [draw, shape=circle,fill=black,scale=0.5] (b3) at  (2,0) {};
\node [draw, shape=circle,fill=black,scale=0.5] (c) at  (4,2) {};
\node [draw, shape=circle,fill=black,scale=0.5] (d) at  (6,2) {};
\node [draw, shape=circle,fill=black,scale=0.5] (e1) at  (8,4) {};
\node [draw, shape=circle,fill=black,scale=0.5] (e2) at  (8,2) {};
\node [draw, shape=circle,fill=black,scale=0.5] (e3) at  (8,0) {};
\node [draw, shape=circle,fill=black,scale=0.5] (f1) at  (10,4) {};
\node [draw, shape=circle,fill=black,scale=0.5] (f2) at  (10,2) {};
\node [draw, shape=circle,fill=black,scale=0.5] (f3) at  (10,0) {};
\draw (a1)--(b1)--(c)--(d)--(e1)--(f1);
\draw (a2)--(b2)--(c)--(d)--(e2)--(f2);
\draw (a3)--(b3)--(c)--(d)--(e3)--(f3);
\node at (0,4.75)  {$1$} {};
\node at (0,2.75)  {$1$} {};
\node at (0,.75)  {$1$} {};
\node at (2,4.75)  {$0$} {};
\node at (2,2.75)  {$0$} {};
\node at (2,.75)  {$0$} {};
\node at (4,2.75)  {$2$} {};
\node at (6,2.75)  {$2$} {};
\node at (8,4.75)  {$0$} {};
\node at (8,2.75)  {$0$} {};
\node at (8,.75)  {$0$} {};
\node at (10,4.75)  {$1$} {};
\node at (10,2.75)  {$1$} {};
\node at (10,.75)  {$1$} {};
\end{tikzpicture}
\caption{A graph for which the lower bound (Lemma~\ref{lowB0}) is attained, with $p=3$. } \label{bistar}
\end{figure}

Notice that this lower bound is sharp, as seen in the graph of Figure~\ref{bistar}.

\begin{coro} \label{lowB0co} Let $G$ be a graph with order $n$ and maximum degree $\Delta \ge 4$. Let $p$ be 
a positive integer such that $3\le p \le \Delta -1$. Let $f=(B_0,B_1,B_2)$ be any $p$-StRD function on $G$.  Then
$$ |B_0| \ge \frac{p}{p-1} \left(n-\gamma_{StR}^p(G)  \vphantom{\frac{}{}}  \right)$$
\end{coro}

Now, we prove the tightness of the upper bound provided by Corollary~\ref{regularesco}. 
To do that, it is sufficient to consider the $(4,5)$-cage graph, known as the Robertson graph. It is 
a $4$-regular graph with $n=19$, and girth $g=5$. Since $\Delta=4$ then $p$ must be $3$. The next 
example shows that $\gamma_{StR}^3(G)=n-p^2+1=11$.

\begin{example}
Let $G$ be the $(4,5)$-cage, the Robertson graph. For this graph, it can be shown that
$ \gamma_{StR}^3(G)=11.$
\end{example}
\noindent {\bf Proof.}
It is readily to prove that $\gamma_{3StR}(G) \leq 11$ by following the construction described in the 
proof of Proposition~\ref{regulares}. 

To see that $\gamma_{StR}^3(G) \ge 11$, we reasoning by contradiction. Assume that $\gamma_{StR}^3(G) \leq 10$. 
Let $f=(V_0,V_1,V_2,V_3)$ be a $\gamma_{StR}^3(G)$-function such that $V_1$ has maximum cardinality. By 
Proposition~\ref{lowbound}, we have that $\gamma_{GStR}^3(G) \geq 7$. Therefore, 
$\gamma_{StR}^3(G)\in \{7,8,9,10\}$. Since $n=19$ and $p=3$, by Corollary~\ref{lowB0co}, we deduce that
$$ |V_0| \ge \left\lceil \frac{3}{2} \left(19- \gamma_{StR}^3(G) \vphantom{\frac{}{}}  \right)  \right\rceil.$$
Clearly, $|V_1| + |V_2| + |V_3| = 19 - |V_0|.$  If $7\le \gamma_{StR}^3(G) \le 8$ then $17\le |V_0| \le 18$ 
which implies $|V_1| + |V_2| + |V_3| \le 2$ and hence $\gamma_{StR}^3(G)\le 6, $ a contradiction. Therefore, 
$9\le \gamma_{StR}^3(G) \le 10.$

Since $f$ is a $\gamma_{StR}^3(G)$-function such that $V_1$ has maximum cardinality, the only 
possibilities are either $\gamma_{StR}^3(G)=10$ with $|V_1|=|V_3|=2,|V_2|=1,|V_0|=14$ or either
$\gamma_{StR}^3(G)=9$ with $|V_1|=|V_2|=1,|V_3|=2,|V_0|=15$. 

Clearly, every vertex with a label $0$ must have a strong neighbor because $f$ is a $3$-StRDF.
Besides, each vertex with a label $3$ is adjacent to, at most, $4$ vertices labeled with $0$ 
and each vertex with a label $2$ is adjacent to, at most, $3$ vertices labeled with $0$.
If $|V_1|=|V_3|=2,|V_2|=1,|V_0|=14$ then we have that
$n=19\le |N[V_3]|+|N[V_2]|+|V_1| \le (2+8)+(1+3)+2=16$, a contradiction. In other case, if
$|V_1|=|V_2|=1,|V_3|=2,|V_0|=15$ then $n=19\le |N[V_3]|+|N[V_2]|+|V_1| \le (2+8)+(1+3)+1=15$, 
again a contradiction. 

$\hfill \Box$

Next, we present a result with a probabilistic approach providing an upper bound. It is described 
in terms of the order, the maximum and minimum degree of the graph and the value of $p$. 


\begin{pro} \label{cota probab}
Let $G$ be a graph with order $n$, minimum degree $\delta $ and maximum degree $\Delta \ge 4$. 
Let $p$ be a positive integer such that $3\le p \le \Delta -1$, such that 
$\left\lceil\frac{\Delta}{p}\right\rceil<\delta$. Then,
$$\gamma_{StR}^p(G) \le \frac{\left( 1+\left\lceil\frac{\Delta}{p}\right\rceil \right) n }{1+\delta} 
\left( \ln \left( \frac{1+\delta}{1+\left\lceil\frac{\Delta}{p}\right\rceil} \right) +1  \right).$$
\end{pro}

\noindent {\bf Proof.}
Let $A\subseteq V(G)$ be a subset of vertices of $G$ and let $\xi\in(0,1)$ be the probability that a 
vertex $v\in V(G)$ belongs to the set $A$.  We assume that two vertices can independently belong to the 
set $A$. Let $B\subseteq V(G)$ be the subset of vertices of $G$ such that do not belong to set $A$ neither 
have neighbors in $A$, that is $B=V(G)-N[A]=(N[A])^c=A^c \cap N(A)^c$. Then, for each vertex $v\in V(G)$ 
we have that
$$
    \begin{array}{rcl}
        P[v\in B] & = & (1-\xi)(1-\xi)^{d(v)} \\ [0.75em]
                  & = & (1-\xi)^{1+d(v)}\leq (1-\xi)^{1+\delta (G)},
    \end{array}
$$
since $0<\xi<1$ and $\delta (G)\le d(v)$, for any $v\in V(G)$. 

Now, for each vertex $v\in V(G)$, we define the following random variable
$$
    X(v)=\left\{
        \begin{array}{ll}
            1+\left\lceil\frac{\Delta}{p}\right\rceil & \mbox{si $v\in A$,} \\[0.5em]
            0 & \mbox{si $v\in N(A)-A$,} \\[0.5em]
            1 & \mbox{si $v\in B=V(G)-N[A]$.}
        \end{array}
          \right. 
$$

It is not difficult to upper bound its expected value, for any $v\in V(G)$, as follows
$$\begin{array}{rcl}
E[X(v)] & = & (1+\left\lceil\frac{\Delta}{p}\right\rceil)P[v\in A]+P[v\in B] \\%
[0.75em]
& = & (1+\left\lceil\frac{\Delta}{p}\right\rceil)\xi+P[v\in B] \\[0.75em]
& \leq & (1+\left\lceil\frac{\Delta}{p}\right\rceil)\xi+(1-\xi)^{1+\delta (G)}%
\end{array}%
$$

\noindent Then, the value that $X(v)$ assigns to each vertex $v\in V(G)$, leads us to a 
function $f:V(G)\rightarrow \{0,1,\ldots ,1+\left\lceil \frac{\Delta}{p}\right\rceil\}$,
such that $f(v)=X(v)$ for each $v\in V(G)$. Since for every vertex $w\in V(G)$, with 
$f(w)=0,$ it has at least one neighbor $u$ in $A$ such that $f(u)=1+\left\lceil
\frac{\Delta}{p}\right\rceil \ge 1+ \left\lceil\frac{1}{p} |N(u)\cap B_0| \right\rceil$
the $f$ is a $p$-StRD function. As a consequence, we have that
$$ 
    \begin{array}{l}
        E[f(V)]  =  \displaystyle\sum_{v\in V(G)}E[f(v)]=\sum_{v\in V(G)}E[X(v)] \\[1.25em]
        \hskip 1.4truecm \leq  \displaystyle\sum_{v\in V(G)}\left((1+\left\lceil\frac{\Delta}{p}\right\rceil)\xi+(1-\xi)^{1+\delta (G)}\right) \\[1.25em]
        \hskip 1.4truecm  =  (1+\left\lceil\frac{\Delta}{p}\right\rceil)n\xi+n(1-\xi)^{1+\delta (G)}
    \end{array}
$$

\noindent Since $0<\xi<1,$ it follows that $(1-\xi)<\mathrm{e}^{-\xi}$ and then
\begin{equation}\label{cotaesperanzaf(V)}  
E[f(V)]\leq \left( 1+\left\lceil\frac{\Delta}{p}\right\rceil\right) n\xi+
n\mathrm{e}^{-\xi(1+\delta (G))}  
\end{equation}

\noindent For each value $\xi\in (0,1)$ minimizing the value of the expression~(\ref{cotaesperanzaf(V)})
it must be
$$\left(1+\left\lceil\frac{\Delta}{p}\right\rceil \right)n-n(1+\delta (G))\mathrm{e}^{-\xi(1+\delta (G))}=0.$$
Therefore, $\mathrm{e}^{-\xi(1+\delta (G))}=\frac{1+\left\lceil \frac{\Delta}{p}\right\rceil}{1+\delta (G)}$
and we deduce that $\xi=\frac{1}{1+\delta(G)}\ln \left( \frac{1+\delta (G)}{1+\left\lceil
\frac{\Delta}{p}\right\rceil}\right).$ 


It is readily to see that $\xi < 1$ because $\ln\left(\frac{1+\delta(G)}{1+\left\lceil
\frac{\Delta}{p}\right\rceil}\right)< \ln \left(\frac{1+\delta (G)}{2}\right)<1+ \delta(G),$
for any $\delta(G).$ Observe also that $\xi>0$ since $\left\lceil\frac{\Delta}{p}\right\rceil <\delta(G)$.
Finally, since $n( 1+ \delta (G))^2 \mathrm{e}^{-\xi(1+ \delta (G))}>0,$ 
we may derive that the critical value of $\xi$ is a local minimum.

\noindent Hence, by using~(\ref{cotaesperanzaf(V)}), we obtain
$$ 
\begin{array}{rl} 
\gamma_{StR}^p(G) \le & \left(1+
\left\lceil\frac{\Delta}{p}\right\rceil \right)\frac{n}{1+ \delta
}\ln\left(\frac{1+ \delta }{1+
\left\lceil\frac{\Delta}{p}\right\rceil} \right)\\
&+\left(1+
\left\lceil\frac{\Delta}{p}\right\rceil \right)\frac{n}{1+
\delta},
\end{array}$$
which concludes the proof.

\hfill $ \Box $

Let us conclude this section with a lower bound expressed in terms of $p$ and the order of the 
graph and a direct consequence for graphs containing a universal vertex.

\begin{pro} \label{lowbound} Let $G$ be a connected graph with order $n$ and maximum degree $\Delta \ge 4$. Let $p$ be a positive integer such that $3\le p \le \Delta -1$. Then
$$ \gamma_{StR}^p(G) \ge  \left\lceil \frac{n+p-1}{p} \right\rceil .$$
If $n \equiv 1 \ (\!\!\!\mod p)$ then equality holds if and only if $\Delta=n-1.$
\end{pro}

\noindent  {\bf Proof.} Let $f=(B_0,B_1,B_2)$ be a $\gamma_{StR}^p(G)$-funcion. Let us denote by $B_0^1$ the 
set of vertices in $B_0$ that have, at most, $p-1$ neighbors in $B_2$ and $B_0^2=B_0-B_0^1.$ 
Clearly, $n=|B_0^1|+|B_0^2|+|B_1|+|B_2|$. Observe that each vertex in $B_1\cup B_2$ contributes with 
one unit, by itself, to the weight of $f$ and each vertex $v \in B_0$ contributes with 
$\frac{|N(v)\cap B_2|}{p}$ to the total weight of $f$. Hence
$$  
    \begin{array}{rcl}
        \gamma_{StR}^p(G) & \ge
        & |B_1| + |B_2| + \displaystyle\sum_{v \in B^1_0} \frac{|N(v)\cap B_2|}{p}\\
        & &+ 
         \displaystyle \sum_{v \in B^2_0} \frac{|N(v)\cap B_2|}{p}    \\[1.25em]
	& \ge &  |B_1| + |B_2| + \displaystyle\frac{1}{p} |B^1_0| + |B_0^2|\\[1.25em] 
    &=&  n - |B^1_0|+ \frac{1}{p} |B^1_0| \\[1.5em]
	& =  &  \displaystyle n- \left(1-\frac{1}{p}\right) |B_0^1| \\[1.25em] 
    &\ge & n-\frac{p-1}{p}(n-1)  = \displaystyle\frac{n+p-1}{p},\\
    \end{array}
$$
since $|B_0^1| \le |B_0| \le n-1$ and $\gamma_{StR}^p(G)$ is an integer.

\noindent Now, let us assume that $n \equiv 1 \ (\!\!\!\mod p)$. On the one hand, if 
$\gamma_{StR}^p(G)= \frac{n+p-1}{p}$, then all previous inequalities became equalities and therefore 
$|B^1_0|=n-1$ and $|B_2|=1$, which implies that $\Delta=n-1$. On the other hand, if 
$\Delta=n-1$ we know that $\gamma_{StR}^p(G) \ge \left\lceil\frac{n+p-1}{p}\right\rceil$. 
To see the other inequality we 
define the function $f$ such that $f(u)=\left\lceil \frac{n-1}{p}  \right\rceil + 1$, for a vertex $u$ such 
that $d_G(u)=n-1$, and $f(v)=0$ for all $v\in N(u)$. Then, $ \gamma_{StR}^p(G) \le w(f)=
\left\lceil \frac{n-1}{p}  \right\rceil +1= \left\lceil \frac{n+p-1}{p} \right\rceil .$

\hfill $\Box$

\begin{coro} \label{universal} Let $G$ be a connected graph with order $n$ and maximum degree $\Delta \ge 4$. Let $p$ be a positive integer such that $3\le p \le \Delta -1$. If $\Delta = n-1$, then
$$ \gamma_{StR}^p(G)  = \left\lceil \frac{n+p-1}{p} \right\rceil = n - \left\lfloor \frac{p-1}{p} \Delta \right\rfloor.$$
\end{coro}


\section{Exact values}


This section is devoted to studying the exact value of the $p$-strong Roman domination number 
in certain families of graphs of interest. We start with the complete bipartite graphs.

\begin{pro}\label{prop:bipartites} Let $p$ be a positive integer such that $3\le p \le \Delta -1$. Let $2\le r \le s$ be 
two positive integers such that $s\ge 4$. Then
$$ \gamma_{StR}^p(K_{r,s})= \left\{
    \begin{array}{ccc}
         2+ \left\lceil \frac{s}{p}  \right\rceil & \text{ if $r=2,$}  \\[1.25 em]
          \left\lceil \frac{r+p-1}{p} \right\rceil + \left\lceil\frac{s+p-1}{p}\right\rceil &
            \text{ if $r\ge 3.$}
    \end{array} 
    \right.
$$
\end{pro}

\noindent {\bf Proof.} Let us denote by $n=n(K_{2,s})=s+2.$ 
To begin with, let us assume that $r=2$.  By applying Proposition~\ref{up1} we
have that $\gamma_{StR}^p(K_{2,s})\le n-\Delta(K_{2,s})+\left\lceil\frac{\Delta(K_{2,s})}{p}
\right\rceil=s+2-s+\left\lceil\frac{s}{p}\right\rceil=\left\lceil\frac{s+2p}{p}\right\rceil.$
If $\gamma_{StR}^p(K_{2,s})\le \left\lceil\frac{s}{p}\right\rceil +1$ then there must be a $p$-StRD 
function $f$ having weight $w(f)\le \left\lceil\frac{s}{p}\right\rceil +1.$
Then, by Corollary~\ref{lowB0co}, the number of vertices labeled with a $0$ can be bounded as follows
$$ 
\begin{array}{rcl}
    |B_0| &\ge& \displaystyle\frac{p}{p-1} \left( n-\left\lceil\frac{s}{p}\right\rceil-1 \right)\\[1.25 em]
    &=&
    \displaystyle\frac{p}{p-1} \left( \left\lfloor \frac{p-1}{p}s\right\rfloor+1 \right)\\[1.25 em] &\ge& \displaystyle\frac{p}{p-1} \ \frac{p-1}{p} s = n - 2.
\end{array}    
$$
Hence, $|B_1|+|B_2|\le 2,$ with $B_2\neq\emptyset$ because $B_0$ is non-empty. We have to consider
several situations.

{\bf Case 1.} If $|B_1|=0$ and $|B_2|=1$ then $|B_0|=n-1$ and $\Delta(K_{2,s})=n-1$, a contradiction.

{\bf Case 2.} If $|B_1|=|B_2|=1$ then $|B_0|=s$ and $N(B_2)=B_0$ which implies that
$w(f)\ge f(B_1)+f(B_2)=1+1+\left\lceil\frac{s}{p}\right\rceil=2+\left\lceil\frac{s}{p}\right\rceil,$
again a contradiction.

{\bf Case 3.} Assume that $|B_1|=0,|B_2|=2$ and let us denote by $B_2=\{u,v\}.$ If $u,v$ are adjacent
then $f(u)+f(v) \ge 1+\left\lceil\frac{s-1}{p}\right\rceil + 2$ implying that
$w(f)\ge 3+\left\lceil\frac{s-1}{p}\right\rceil \ge 2+\left\lceil\frac{s}{p}\right\rceil$ which is
not possible. If $u,v$ are not adjacent then $|N(\{u,v\})|=s$ and $B_0=N(u)=N(v)$, because $f$ is 
a $p$-StRDF. Hence, $f(u)+f(v)\ge 2 \left( 1+\left\lceil\frac{s}{p}\right\rceil \right)>
2+\left\lceil\frac{s}{p}\right\rceil,$ which is a contradiction.

So, it must be $\gamma_{StR}^p(K_{2,s})\ge 2+\left\lceil\frac{s}{p}\right\rceil$ and the equality is
proven.

Now, let us suppose that $r\ge 3.$ Let $u$ be a vertex belonging to the $r$-class and $v$ be a vertex
of the $s$-class. The function defined as $f(u)=1+\left\lceil\frac{s-1}{p}\right\rceil,
f(v)=1+ \left\lceil\frac{r-1}{p}\right\rceil $ is a $p$-StRD function in $K_{r,s}.$ Therefore
$\gamma_{StR}^p(K_{r,s})\le  \left\lceil\frac{r+p-1}{p}\right\rceil+ \left\lceil\frac{s+p-1}{p}\right\rceil. $

Reasoning by contradiction, let us assume that
$\gamma_{StR}^p(K_{r,s}) \le  \left\lceil\frac{r-1}{p}\right\rceil+ \left\lceil\frac{s-1}{p}\right\rceil +1. $
Let $f=(B_0,B_1,B_2)$ be a $\gamma_{StR}^p$-function. Again by Corollary~\ref{lowB0co}, we have that
$$  
\begin{array}{rl}
    |B_0|  \ge & 
     \frac{p}{p-1}\left(r+s-\left\lceil\frac{r-1}{p}\right\rceil - \left\lceil\frac{s-1}{p}\right\rceil -1 \right)
    \\[1.25em]
     =& \frac{p}{p-1}\left(\left\lfloor\frac{p-1}{p}(s-1)\right\rfloor \right.+ 1 \\
     &\left.+\left\lfloor\frac{p-1}{p}(r-1)\right\rfloor 
    +1-1\right) \\[1.25em]
     \ge& \frac{p}{p-1}\left(\frac{p-1}{p}(s-1) +\frac{p-1}{p}(r-1)-1\right) \\[1.25em]
     = &  r+s-2-\frac{p}{p-1}
    \end{array}
$$
As the function $h(p)=\frac{p}{p-1}$ is a non-increasing function for positive values of $p$, and since $h(3)=\frac{3}{2}$ we 
deduce that $|B_0|\ge r+s-\frac{7}{2},$ which in turn lead us to $|B_1|+|B_2|\le 3.$ Since $B_2\neq \emptyset$ then we have 
to consider different cases: $|B_1|=j$ and $1\le |B_2|\le 3-j$, for all $j\in\{0,1,2\}$.

{\bf Case 1.} If $|B_1|=2,|B_2|=1$ then there must be $r=3, N(B_2)=B_0$ and the $3$-class of $K_{3,s}$ coincides with
$B_1\cup B_2$. But in this case, it would be
$$  w(f) =2+1+\left\lceil\frac{s}{p}\right\rceil > \left\lceil\frac{r-1}{p}\right\rceil+ \left\lceil\frac{s-1}{p}\right\rceil +1  $$

{\bf Case 2.} $|B_1|=1$ and $1\le |B_2|\le 2$. If $|B_1|=|B_2|=1$ then $|B_0|=n-2$ which is impossible because $r\ge 3$ 
and $f$ is an StRDF. Hence, $|B_1|=1$ and $|B_2|=2.$ Let us denote by $B_1=\{z\}, B_2=\{u,v\}.$ If $u,v$ are not adjacent 
then it must be $r=3$ and $B_1\cup B_2$ is the $3$-class of $K_{3,s},$ because $f$ is an StRDF.  Then, we deduce that
$w(f)= 2\left( 1+\left\lceil\frac{s}{p}\right\rceil \right)+1 > 
\left\lceil\frac{r-1}{p}\right\rceil+ \left\lceil\frac{s-1}{p}\right\rceil +1,$ a contradiction. If $u,v$ are adjacent then, 
without loss of generality, we may assume that $d(z)=s.$ Hence,
$w(f)= f(z)+f(u)+f(v)=1+\left\lceil\frac{s-1}{p}\right\rceil+1+\left\lceil\frac{r-2}{p}\right\rceil+1 > 
\left\lceil\frac{r-1}{p}\right\rceil+ \left\lceil\frac{s-1}{p}\right\rceil +1.$

{\bf Case 3.} $|B_1|=0$ and $1\le |B_2|\le 3$. As $r\ge 3$ then $\Delta \le n-3$ and therefore $B_2\ge 2.$ If the induced 
subgraph by the vertices of $B_2$ is an edgeless subgraph then $r=|B_2|=3$ and $w(f)=
3\left( 1+\left\lceil\frac{s}{p}\right\rceil \right)+1>\left\lceil\frac{r-1}{p}\right\rceil+ \left\lceil\frac{s-1}{p}\right\rceil+1.$ 
Then, there must be adjacent vertices in the set $B_2.$ If $|B_2|=2$ then $w(f)=f(B_2)\ge 
1+\left\lceil\frac{s-1}{p}\right\rceil+1+\left\lceil\frac{r-1}{p}\right\rceil$, a contradiction. If $|B_2|=3$ then
$ w(f) = f(B_2) \ge 
1+\left\lceil\frac{s-1}{p}\right\rceil+1+\left\lceil\frac{r-2}{p}\right\rceil+1 =
2+\left\lceil\frac{s-1}{p}\right\rceil+\left\lceil\frac{r+p-2}{p}\right\rceil \ge
2+\left\lceil\frac{s-1}{p}\right\rceil+\left\lceil\frac{r-1}{p}\right\rceil $, again a contradiction. So the result
holds.

\hfill $\Box$

Our next result provides the exact value of the $p$-strong Roman domination number for bi-star graphs. The proof 
is quite similar to that of Proposition~\ref{prop:bipartites}, so we leave the details to the reader.

\begin{pro} Let $r, s$ be two integers such that $1 \le r \le s$. Let $T_{r,s}$ be the bi-star 
graph with order $n=r+s+2$ and maximum degree $\Delta \ge 4$. Let $p$ be a positive integer such 
that $3\le p \le \Delta -1$. Then,
$$ \gamma_{StR}^p(T_{r,s})= 2+ \left\lceil \frac{r}{p}  \right\rceil +  \left\lceil \frac{s}{p}  \right\rceil.$$
\end{pro}

We conclude by characterizing those graphs having the smallest possible values of the 
$p$-strong Roman domination number. 

\begin{pro} Let $G$ be a graph with order $n$ and maximum degree $\Delta \ge 4$. Let $p$ be a positive 
integer such that $3\le p \le \Delta -1$. Then $\gamma_{StR}^p(G)=3$ if and only if $G=K_1 \vee H$, where 
$p+2\le n\le 2p+1$, $\Delta=n-1$ and $H$ is any graph with $n-1$ vertices.
\end{pro}

\noindent {\bf Proof.}
Assume that $\gamma_{StR}^p(G)=3$. Let $f=(B_0,B_1,B_2)$ be a $p$-StRD function on $G$ with minimum 
weight $w(f)=3$. Since $w(f)=|B_1|+\sum_{v\in B_2}f(x)$, hence there are only two possibilities: (i) $B_1=\emptyset$ 
and  $B_2=\{v\}$, with $f(v)=3$ or (ii) $|B_1|=1$ and  $B_2=\{v\}$, with $f(v)=2$. \\
\noindent (i) If $B_1=\emptyset$ and $B_2=\{v\},$ with $f(v)=3,$ then $B_0=V \smallsetminus \{ v \}$, that is, 
every vertex in $B_0$ must be adjacent to $v$ and $p+1\le |B_0|\le 2p$, since $f(v)=3$. Therefore, $p+2\le n\le 2p+1$ 
and $\Delta=n-1$, and then, $G=\{v\} \vee H,$ where $H$ is a subgraph on $p+1\le |V(H)| \le 2p$ vertices. \\
\noindent (ii) If $B_1=\{u\}$ and $B_2=\{v\},$ with $f(v)=2$, then $B_0=V \smallsetminus \{ u, v \}$ and all vertices 
in $B_0$ must be adjacent to $v$, hence $1\le |B_0|\le p$, since $f(v)=2$, and $3 \le n \le p+2$. Now, we distinguish 
two subcases: $|N(v)|=n-1$ or $|N(v)|=n-2$.\\
Suppose that $|N(v)|=n-1$. If $|B_0| < p$, then, we can label the vertex $u$ with $0$ and we have $\gamma_{StR}^p(G)=2$, 
which is a contradiction. Hence, $|B_0|=p$ and there exists a $p$-StRD function $f$ on $G$ with minimum weight $w(f)=3$, 
where $f(v)=3$, which lead us to case (i). \\
Suppose now that $|N(v)|=n-2$. Since $1\le |B_0|\le p$, then $n \le p+2$. Due to $3\le p \le \Delta -1$, then 
$p+1 \le \Delta$. We deduce that $n-1 \le \Delta$ and, therefore, $n-1 = \Delta$, which means that there exists a 
vertex $w \in B_0$ such that $d(w)=\Delta=n-1$, which describes the graph of the case (i). The reciprocal is trivial.
$\hfill \Box$

\begin{pro}  Let $G$ be a graph with order $n$ and maximum degree $\Delta \ge 4$. Let $p$ be a positive integer 
such that $3\le p \le \Delta -1$. Then $\gamma_{StR}^p(G)=4$ if and only if one of the following conditions hold
\begin{enumerate}
\item $\Delta=n-1$, $2p+2 \le n \le 3p+1$ and $G= K_1 \vee H$ where $H$ is any graph on $2p+1\le |V(H)| \le 3p$ vertices.
\item $\Delta=n-2$, $4 \le n \le 2p+2$ and $G= H_1 \vee H_2$, where $H_1 \subseteq K_2$ and $H_2$ is a subgraph 
on $2\le |V(H_2)| \le 2p$ vertices.
\item $\Delta=n-2$, $p+3\le n\le 2p+2$ and $G= K_1 \vee H$, where $K_1=\{z\}$, with $d(z)=\Delta(G)$, and $H$ 
is a subgraph on $|V(H)|=n-1$ vertices.
\end{enumerate}
\end{pro}

\noindent {\bf Proof.}
Assume that $\gamma_{StR}^p(G)=4$. Let $f=(B_0,B_1,B_2)$ be a $p$-StRD function on $G$ with minimum weight 
$w(f)=4$. Since $w(f)=|B_1|+\sum_{v\in B_2}f(x)$, hence there exists different possibilities. \\
\noindent {\bf Case 1:} Assume that $B_1=\emptyset$. 

{\bf Subcases 1a:} Suppose that $B_2=\{v\}$, with $f(v)=4$. Then every vertex in $B_0$ must be adjacent to $v$ and 
$2p+1\le |B_0|\le 3p$, since $f(v)=4$. Therefore, $2p+2\le n\le 3p+1$ and $\Delta=n-1$, and then, $G=\{v\} \vee H$, 
where $H$ is a subgraph on $2p+1\le |V(H)| \le 3p$. 

{\bf Subcase 1b:} Suppose that $B_2=\{v, w\}$, with $f(v)=f(w)=2$. Then every vertex in $B_0$ must be adjacent to $v$ 
and $w$, then $2 \le |B_0|\le 2p$. Therefore, $4 \le n \le 2p+2$ and $\Delta=n-2$, and then, $G= H_1 \vee H_2$, 
where $H_1 \subseteq K_2$ and $H_2$ is a subgraph on $2\le |V(H_2)| \le 2p$ vertices.

\noindent {\bf Case 2:} Assume that $B_1 \neq \emptyset$. 

{\bf Subcase 2a:} Suppose that $B_1=\{u\}$ and $B_2=\{v\},$ with $f(v)=3$. Then every vertex in $B_0$ must be adjacent 
to $v$ and $p+1\le |B_0|\le 2p$, since $f(v)=3$, and $p+3\le n\le 2p+2$. If $w \in N(v)$, then $\Delta=n-1$ and there 
exists a $p$-StRD function $f$ on $G$ with minimum weight $w(f)=4$, where $f(v)=4$, which lead us to case (1a). If 
$w \notin N(v)$, then $|B_0|\le n-2$, since $f(v)=3$, $\Delta=n-2$ and $G= K_1 \vee H$, where $K_1=\{z\}$, with 
$d(z)=\Delta(G)=n-2$, and $H$ is a subgraph on $|V(H)|=n-1$ vertices.

{\bf Subcase 2b:} Suppose that $B_1=\{u, w\}$ and $B_2=\{v\},$ with $f(v)=2$. Then every vertex in $B_0$ must be 
adjacent to $v$ and $1\le |B_0|\le p$, since $f(v)=2$, and $4 \le n\le p+3$. If $u,w \in N(v)$, then $d(v)=n-1=p+2$ 
and $n=p+3$, therefore, $|B_0|=p$, which implies that there exists a $p$-StRD function $f$ on $G$ with minimum weight 
$w(f)=3$, where $f(v)=3$, which is a contradiction, since we are assuming that $\gamma_{StR}^p(G)=4$. If $u,w \notin N(v)$, 
then $1\le |B_0|\le p$ and $\Delta = n-3$, hence, $p \le \Delta -1 = n-4$, that is,  $p+4 \le n$, which is a contradiction, 
since $4 \le n\le p+3$. If $u \in N(v)$ and $w \notin N(v)$, then necessarly $|B_0|=p$, since $f(v)=2$. Therefore $n=p+3$ 
and $\Delta=p+1 = n-2$, which leads us to case (2a).

\noindent The reciprocals are trivial.
$\hfill \Box$

{ 
\section{ Conclusion}
Recently, many definitions of Roman domination models for graphs have been proposed. 
In this work, we introduce the concept of p-strong Roman domination, which enables the 
development of newer, more adaptable, and less expensive defensive strategies.

The NP-completeness of the problem has been explored for bipartite and chordal graphs by 
linking it to the Exact 3-Cover problem. Various general upper and lower bounds have been 
examined, along with an upper bound derived using probabilistic methods. Concerning the 
study of exact values, specific cases like Robertson's (4,5)-cage, where 3-StR equals $11$, 
and extensive families of graphs such as complete bipartite graphs or bi-stars have been investigated.
}

{ 
\section{Future research directions.}
This work opens several compelling avenues for future research, particularly in light of 
the NP-completeness of the p-strong Roman domination problem.

Several promising directions warrant investigation. Firstly, seeking tighter bounds, either 
by improving existing ones or by considering other graph invariants, would be valuable. 
Another interesting possibility is exploring new inequalities that relate this parameter to 
parameters of other domination types beyond those mentioned in Remark~\ref{otras}.

Regarding Remark~\ref{otras}, characterizing the graphs that achieve lower and upper bounds 
would be a significant contribution. Determining the exact value of the p-strong Roman 
domination number in other graph families and graph products is also a worthwhile pursuit.

Furthermore, studying the exact value of the p-strong Roman domination number in trees 
with specific structures is of interest. Based on these findings, an attempt to establish 
a general bound for any tree of order $n\ge 5$ could be undertaken.
}

\section*{Author contributions statement} 
All authors contributed equally to this work.

\section*{Conflicts of interest} 
The authors declare no conflict of interest.

\section*{Data availability} 
No data was used in this investigation.


\end{document}